\documentclass[12pt]{article}
\usepackage[dvips]{graphicx}
\usepackage[cp1251]{inputenc}
\usepackage[russian]{babel}
\usepackage{amsmath,amssymb,amsthm}

\newcommand{\myincludegraphics}[1]{\includegraphics{#1}}

\makeatletter
\def\@seccntformat#1{\csname the#1\endcsname. }
\def\@allowhyphens{\penalty\@M\hskip\z@skip}
\def\set@low@box#1{\setbox\tw@\hbox{,}\setbox\z@\hbox{#1}\dimen\z@\ht\z@
     \advance\dimen\z@ -\ht\tw@
     \setbox\z@\hbox{\lower\dimen\z@ \box\z@}\ht\z@\ht\tw@ \dp\z@\dp\tw@ }
\def\@glqq{{\ifhmode \edef\@SF{\spacefactor\the\spacefactor}\else
     \let\@SF\empty \fi \leavevmode
     \set@low@box{''}\box\z@\kern-.04em\@allowhyphens\@SF\relax}}
\def\glqq{\protect\@glqq}
\def\@grqq{\ifhmode \edef\@SF{\spacefactor\the\spacefactor}\else
     \let\@SF\empty \fi \kern-.07em``\kern.07em\@SF\relax}
\def\grqq{\protect\@grqq}
\makeatother

\newtheorem{thm}{Теорема}[section]

\newtheorem{prop}[thm]{Предложение}
\newtheorem{fig}[thm]{Рисунок}
\newtheorem{Ex}[thm]{Пример}

\newtheorem{lemma}[thm]{Лемма}
\newtheorem{cor}[thm]{Следствие}

\theoremstyle{definition}
\newtheorem{rem}[thm]{Замечание}
\newtheorem{dfn}[thm]{Определение}

\newcommand{\bw}{{Ч$\cap$Б}}

\newcommand{\divs}{\mathrel{.\hbox to 0pt{\hss\raisebox{1.2ex}[0pt][0pt].\hss\raisebox{.6ex}[0pt][0pt].}}}
\newcommand{\uline}[1]{{\hbox to 0pt{\underline{\hphantom{{#1}}}\hss}{#1}}}

\begin{document}
\title{
УДК 514.172.45\\
Некоторые примеры смежностных\\
многогранников коразмерности 4}
\author{Ростислав Девятов}
\date{}
\maketitle

{
\begin{abstract}

В настоящей статье построена серия комбинаторныx типов смежностных многогранников в
$\mathbb R^{2d}$ с числом вершин $N=2d+4$. Все построенные многогранники имеют плоскую диаграмму
Гейла определенного вида, а именно: в ней ровно $d+3$ черных точек, лежащих в выпуклом
положении.

Такие диаграммы Гейла параметризуются 3-деревьями (деревьями с
некоторой дополнительной структурой).

Для всех многогранников  построенной серии число граней
размерности $m$, содержащих вершину $A$, зависит лишь от $d$ и $m$.
\end{abstract}
}

\section{Введение}
\subsection*{Смежностные и циклические многогранники}
Многогранник в $\mathbb R^D$ называется \textit{смежностным}, если любые $\lfloor D/2\rfloor$
его вершин образуют грань. Мы будем рассматривать случай, когда $D$ четно ($D=2d$), тогда
любые $d$ вершин смежностного многогранника образуют  грань.

Примером смежностного многогранника является циклический многогранник \cite{zeig}. Для его
построения используется \textit{кривая моментов в $\mathbb R^D$}, которая содержит все точки вида
$(t, t^2, t^3, \ldots, t^D)$, где $t\in\mathbb R$. Выпуклая оболочка любого конечного подмножества кривой
моментов называется циклическим многогранником. Циклический многогранник является
смежностным.

Известно \cite{zeig}, \cite{mcmul}, что число граней циклического многогранника фиксированной
размерности максимально возможное для многогранников такой размерности с таким количеством
вершин.

Кроме циклического, даже в четной размерности существуют другие примеры смежностных
многогранников с $D+4$ вершинами (в нечетной размерности ситуация проще, см.~\cite{she}).

Первые примеры смежностных многогранников были получены Грюнбаумом~\cite{gru}:
для $D=4$ есть 2 комбинаторных типа смежностных многогранников, не считая циклического.\par

\subsection*{Диаграмма Гейла (см.~\cite{zeig})}

Векторная диаграмма Гейла позволяет конфигурации $X^*$ из $N$ точек в $D$-мерном пространстве
поставить в соответствие конфигурацию $\overline X$ из $N$ ненулевых векторов в $(N-D-1)$-мерном
пространстве.

При этом комбинаторика $\overline X$ 
полностью определяется комбинаторикой исходной конфигурации
$X^*$ 
(см. определения в \cite{zeig}) 
и наоборот. Это каноническое соответствие, т.~е. 
автоморфизмы $\overline X$ индуцируют автоморфизмы $X^*$ и наоборот.

Центральная проекция (из начала координат) конфигурации $\overline X$ на произвольно выбранную
гиперплоскость $e$ (см.~рис.~\ref{vecaffdiag})
\begin{fig}
$$
\myincludegraphics{szsfigures.1}\label{vecaffdiag}
$$
\end{fig}

\noindent дает конфигурацию  $X$ черных и белых точек в $(N-D-2)$-мерном пространстве~--- аффинную
диаграмму Гейла конфигурации $X^*$.

Аффинная диаграмма Гейла не определяется однозначно конфигурацией $X^*$, т.~к. имеется свобода
выбора гиперплоскости $e$. Например,
две аффинные диаграммы Гейла (см.~рис.~\ref{2aff}) соответствуют одному и тому же комбинаторному типу
исходной конфигурации $X^*$.

\begin{fig}
$$
\begin{array}{cc}
\myincludegraphics{szsfigures.2} & \hspace{5em}\myincludegraphics{szsfigures.3}\label{2aff}
\end{array}
$$
\end{fig}

Если не оговорено противное, в дальнейшем  под словами \glqq{}диаграмма Гейла\grqq{} будет пониматься
именно аффинная диаграмма.

\begin{dfn}
Конечный набор черных и белых точек на плоскости, аффинно
порождающий плоскость, называется (плоской) \textit{диаграммой}.
\end{dfn}

\begin{dfn}
Две плоских диаграммы,  как правило (см. ~\cite{zeig}),
называются \textit{комбинаторно эквивалентными}, если между этими двумя множествами существует биекция,
сохраняющая цвет точек, ориентации троек и отношение \glqq лежать между\grqq{} для троек точек, 
лежащих на одной прямой.\label{combequivzeig}
\end{dfn}

\begin{dfn}
Однако в данной статье мы называем диаграммы \textit{комбинаторно эквивалентными}, если они
комбинаторно эквивалентны в смысле определения \ref{combequivzeig} с точностью до зеркальной 
симметрии.\label{combequivdef}
\end{dfn}

Аналогично, комбинаторную эквивалентность векторных диаграмм (наборов векторов в $\mathbb R^3$) 
мы также будем понимать с точностью до смены ориентации всего пространства.

\begin{dfn}
Пусть  $X$~--- плоская диаграмма, точки которой лежат в общем положении. Будем говорить, что
$M \subseteq X$ обладает свойством \bw, если относительные внутренности множеств
$$
Conv(\{ A \in M: \text{$A$~--- черная точка}\})
$$
и
$$
Conv(\{ A \in M: \text{$A$~--- белая точка}\})
$$
имеют непустое пересечение. (Здесь $Conv (\cdot )$ обозначает выпуклую оболочку.)
\end{dfn}

Нам потребуются следующие известные факты о диаграммах Гейла (см. \cite{zeig}, \cite{she}):

\begin{enumerate}
\item Имеется каноническая биекция между точками диаграмм $X^*$ и $X$. Точку, соответствующую
точке $A  \in X$, будем обозначать через $A^*$.

\item Четномерный смежностный многогранник симплициален, а следовательно, точки
диаграммы Гейла множества вершин смежностного многогранника (для краткости в дальнейшем
будем писать \glqq точки диаграммы Гейла смежностного многогранника\grqq) лежат в общем положении.

\item Точки $X^*$ лежат в общем положении тогда и только тогда, когда точки $X$ лежат в
общем положении.

\item Плоская диаграмма $X$ является диаграммой Гейла некоторого выпуклого многогранника
(т.~е. точки  конфигурации $X^*$ лежат в выпуклом положении)
тогда и только тогда, когда для любой точки $A \in X$ свойство \bw{} выполнено для
множества $X \setminus \{A\}$.

\item Критерий грани многогранника на языке диаграммы Гейла.

Пусть точки конфигурации  $X^*$ являются вершинами некоторого выпуклого многогранника $K$.
Множество $M^*$
является множеством вершин некоторой грани $K$ тогда и только тогда,
когда для множества $X\setminus M$ выполнено свойство \bw.

В связи с этим в дальнейшем иногда будет использоваться следующая терминология: мы 
будем говорить об \glqq удалении\grqq{} некоторых точек из диаграммы, а затем для \glqq 
оставшихся\grqq{} точек проверять свойство \bw. Оно будет выполнено тогда и только тогда,
когда \glqq удаленные\grqq{} точки образуют грань.

\item Комбинаторный тип диаграммы Гейла циклического многогранника размерности $D=2d$ с $D+4$
вершинами таков (см.~рис.~\ref{cyclicfig}): $D+4$ точки лежат на окружности, причем эти точки попеременно
белые и черные.

Следовательно, группа комбинаторных автоморфизмов циклического многогранника действует на множестве
его вершин транзитивно.

\begin{fig}
$$
\myincludegraphics{szsfigures.4}\label{cyclicfig}
$$
\end{fig}

\item Смежностные многогранники \textit{комбинаторно жесткие}, т.~е. комбинаторный тип
смежностного многогранника определяет комбинаторику конфигурации вершин. Следовательно,
комбинаторный тип смежностного многогранника однозначно определяет комбинаторику
его векторной диаграммы Гейла.

\end{enumerate}

\subsection*{Циклический порядок}
Рассмотрим окружность $S$ с расставленными на ней через равные расстояния $n$ точками и с заданной
ориентацией (т.~е. одно из двух направлений ее обхода названо направлением \glqq по часовой
стрелке\grqq, а другое --- \glqq против часовой стрелки\grqq).

\begin{dfn}
Задать на множестве $M$ из $n$ элементов \textit{циклический порядок} значит поставить его
элементы во взаимно однозначное соответствие с выбранными точками на окружности, с точностью до
поворота окружности на угол $2k\pi / n$.
\end{dfn}

Про элементы множества, на котором задан циклический порядок можно говорить о соседних с ними
элементах.

\textit{Противоположным} к данному циклическому порядку называется циклический порядок, в котором
направления по и против часовой стрелки поменялись местами. Этот порядок получится из исходного,
если у $S$ изменить ориентацию.

\begin{dfn}
Задать на множестве $M$ из $n$ элементов \textit{неориентированный циклический порядок} значит
задать на нем циклический порядок с точностью до перехода к противоположному.
\end{dfn}

Циклический порядок на множестве порождает естественным образом неориентированный
циклический порядок.
\subsection*{Основные результаты}
Мы строим смежностные многогранники с $D+4$ вершинами, а следовательно, с двумерной
диаграммой Гейла.

\begin{dfn} Назовем плоскую диаграмму $X$ из $D+4$ черных и белых точек \textit{Т-диаграммой}
(см. примеры Т-диаграмм на рис.~\ref{FirstSamples}), если
\begin{enumerate}
\item $X$ содержит ровно $d+3$ черных точки, причем они лежат в выпуклом положении.
На множестве черных точек естественным образом задан неориентированный циклический
порядок, порожденный обходом границы $(d+3)$-угольника.

\item $X$ является диаграммой Гейла смежностного многогранника.
\end{enumerate}
Смежностный многогранник, соответствующий Т-диаграмме, назовем \textit{Т-многогранником}.
\end{dfn}

\begin{dfn}
Дерево с дополнительной структурой назовем \textit{3-деревом}, если
выполнены следующие условия:

\begin{enumerate}
\item Все его вершины имеют степень либо 1, либо 3.
\item На каждых 3 ребрах, выходящих из вершины степени 3, задан циклический порядок.
Этот порядок естественным образом индуцирует циклический порядок на вершинах степени 1.
\end{enumerate}
\end{dfn}
(См. примеры 3-деревьев на рис.~\ref{treesample})

\begin{dfn}
3-дерево называется \textit{характеристическим деревом} данной Т-диаграммы, если
можно установить взаимно-однозначное соответствие между его вершинами и точками
Т-диаграммы так, что:
\end{dfn}
\begin{enumerate}
\item Белые точки диаграммы соответствуют вершинам степени 3 дерева, а черные~---
вершинам степени 1.
\item Неориентированный циклический порядок на вершинах степени 1 в дереве совпадает с
неориентированным циклическим порядком на черных точках, задаваемым обходом границы
$(d+3)$-угольника.
\item Белая точка $B$ на диаграмме лежит внутри треугольника, образованного тремя
черными точками $A_i, A_j, A_k$, тогда и только тогда, когда из вершины дерева, соответствующей
$B$, к 3 вершинам, соответствующим $A_i, A_j, A_k$, ведут 3 непересекающихся пути.
\newcounter{r117}\setcounter{r117}{2}\refstepcounter{r117}\label{triangs}
\end{enumerate}

\begin{thm}(См. теорему \ref{characteristictree} и лемму \ref{combequiv}.)
\begin{enumerate}
\item Для любой Т-диаграммы существует характеристическое дерево. Оно единственно
с точностью до одновременной смены циклического порядка во всех вершинах на противоположный порядок.
\item Для любого 3-дерева с числом вершин $D+4$, где $D\ge 4$, существует Т-диаграмма, для которой
данное дерево является характеристическим. Она единственна с точностью до комбинаторной
эквивалентности.\qed
\end{enumerate}
\end{thm}

\begin{thm} (См. теорему \ref{equivTpolyTdiag} и лемму \ref{combequiv}.)
Комбинаторный тип Т-многогранника однозначно определяет комбинаторный тип его Т-диаграммы Гейла.\qed
\end{thm}

\begin{cor} (См. следствие \ref{TpolyTree})
\begin{enumerate}
\item Любому Т-многограннику соответствует 3-дерево (характеристическое дерево для его Т-диаграммы
Гейла). Оно единственно с точностью до смены циклического порядка в вершинах на
противоположный.
\item Любое 3-дерево с не менее, чем 8 вершинами, задает единственный с точностью до комбинаторной
эквивалентности Т-многогранник.\qed
\end{enumerate}
\end{cor}

\begin{thm} (См. теорему \ref{EqualIncedence})
В любом Т-многограннике размерности $D$ число граней с
$m$ вершинами, содержащих данную вершину, зависит только от $D$ и $m$,
но не от многогранника и выбранной вершины. Оно совпадает с соответствующим числом
для циклического многогранника.\qed
\end{thm}

\begin{dfn}
Назовем ребро Т-многогранника \textit{замечательным}, если любой набор из $d+1$ вершины,
содержащий его концы, порождает грань.
\end{dfn}

\begin{prop} (См. предложение ~\ref{GoodCrit} и замечание ~\ref{WcorrBintree})
Замечательные ребра смежностного многогранника соответствуют ребрам, инцидентным висячим
вершинам соответствующего 3-дерева.\qed
\end{prop}

\subsection*{Благодарности}
Автор благодарен организаторам летней школы \glqq Современная математика\grqq, где был прочитан
курс о диаграммах Гейла, а также лично лектору этого курса Г.~Ю.~Паниной.\par

\section{Критерий смежностности на языке Т-диаграмм}

\begin{lemma}
В Т-диаграмме все белые точки находятся строго внутри $(d+3)$-угольника, образованного черными
точками.\label{inside}
\end{lemma}
Доказательство. Здесь и далее, мы пользуемся критерием грани для диаграммы Гейла (см. введение).
Пусть белая точка $B$ находится вне $(d+3)$-угольника. Удалим остальные $d$ белых точек.
Множество из оставшихся $d+4$ точек (все черные точки и точка $B$)
не обладает свойством \bw{}, что противоречит
смежностности.\qed

\begin{dfn}
Будем называть \textit{луночкой} Т-диаграммы выпуклую оболочку нескольких подряд идущих вершин
$(d+3)$-угольника.
\end{dfn}

\begin{fig}
$$
\myincludegraphics{szsfigures.5}
$$
\end{fig}

\begin{lemma}
В любом треугольнике, образованном 3 черными точками  Т-диаграммы, лежит ровно одна белая
точка.\label{less2}
\end{lemma}
Доказательство. Докажем вначале, что каждый такой треугольник содержит не более одной белой точки.
Допустим, что в треугольнике $T$  находятся хотя бы 2 белые точки. Удалим из $(d+3)$-угольника
треугольник $T$. При этом $(d+3)$-угольник распадется на три (возможно, вырожденные) луночки.

Пусть у этих луночек соответственно $a$, $b$ и $c$ черных вершин (будем считать, что у вырожденной
луночки 2 вершины). Очевидно, $a+b+c=(d+3)+3$. Внутри луночек в совокупности содержится
не более $d-1$ белой точки.  Поэтому в некоторой из наших трех луночек,
граница которой состоит из $x$ вершин, содержится менее $x-2$ белых точек.
В противном случае во всех луночках в сумме было бы не меньше $(a-2)+(b-2)+(c-2)=d$ белых точек.
Рассмотрим черные точки, которые не лежат на границе выбранной луночки (их $d+3-x$) и
белые точки внутри нее (их не более $x-3$). Мы выбрали не более $d$ точек. (Если мы выбрали
строго меньше, добавим произвольные белые точки так, чтобы всего получилось $d$.)
Указанные $d$ точек не образуют грань, т.~к. свойство \bw{} для оставшихся точек не
выполнено.

\begin{fig}
$$
\myincludegraphics{szsfigures.6}
$$
\end{fig}

Докажем теперь существование белой точки внутри любого треугольника $T$ с черными вершинами.
Разобьем $(d+3)$-угольник на треугольники, так чтобы
среди треугольников разбиения был и треугольник $T$.
Мы получим $d+1$ треугольник, в которых как-то расставлена $d+1$ белая точка.
Утверждение следует из принципа Дирихле.\qed

\begin{prop}
Пусть имеется плоская диаграмма $X$ из $2d+4$ черных и белых точек, обладающая следующими
свойствами:\label{equal1enough}
\begin{enumerate}
\item $D=2d\ge4$\label{ge4};
\item точки $X$ лежат в общем положении;
\item $X$ содержит ровно $d+3$ черных точек, лежащих в выпуклом положении;
\item вне выпуклой оболочки черных точек нет белых точек;
\item внутри каждого  треугольника c черными вершинами лежит ровно одна белая точка.
\end{enumerate}
Тогда $X$~--- Т-диаграмма.
\end{prop}

\begin{fig}
$$
\begin{array}{cc}
\myincludegraphics{szsfigures.7} & \myincludegraphics{szsfigures.8}\label{FirstSamples}\\[10pt]
\myincludegraphics{szsfigures.9} & \myincludegraphics{szsfigures.10}
\end{array}
$$
\end{fig}

Доказательство. Достаточно доказать, что при удалении из $X$ любого набора из $d$ точек, оставшееся
множество обладает свойством \bw{}. Действительно, пусть удалены некоторые
$k$ черных точек и $d-k$ белых точек. Тогда у нас осталось $d-k+3$ черных точек и
$k+1$ белая точка. Разобьем $(d+3)$-угольник на треугольники, так чтобы
$(d-k+3)$-угольник из оставшихся черных вершин целиком состоял из треугольников разбиения.
Всего получился $d+1$ треугольник, причем $d-k+1$ из них составляют наш $(d-k+3)$-угольник.
Мы удалили всего $d-k$ белых вершин, значит внутри какого-то из этих $d-k+1$
треугольников осталась белая точка.\qed

\begin{cor}
Примеры, приведенные на рис.~\ref{FirstSamples}, являются Т-диаграммами Гейла.
\end{cor}

\section{Т-диаграммы и 3-деревья}

\begin{dfn}
Плоская диаграмма из $2d+4$ точек
называется \textsl т\textit{-диаграммой}, если она удовлетворяет условиям 2--5 из предложения \ref{equal1enough}, 
а $d\ge 0$.
\end{dfn}

Пусть имеется т-диаграмма из $2d+4$ точек. Обозначим черные точки в порядке обхода по
часовой стрелке через $A_1, A_2, \ldots, A_{d+3}$. Иногда мы будем называть точку
$A_{d+3}$ точкой $A_0$, а точку $A_1$~--- точкой $A_{d+4}$.

\begin{dfn}
Назовем \textit{граничным треугольником точки $A_i$} треугольник $A_{i-1}A_iA_{i+1}$.
\end{dfn}

Согласно лемме \ref{less2}, в каждом граничном треугольнике содержится ровно одна
белая точка. Про белую точку, которая лежит внутри $\triangle A_{i-1}A_iA_{i+1}$ будем говорить,
что она \textit{соответствует вершине $A_i$.} Черной точке всегда соответствует ровно
одна белая (но не наоборот).

\begin{rem}
Если белая точка соответствует некоторой черной точке $A_i$, то она
лежит внутри только тех треугольников с черными вершинами, у которых
одна из вершин $A_i$.\label{insideonly}
\end{rem}

\begin{lemma}
\textbf{\textup{Лемма об индукции для т-диаграммы.}} Пусть \textsl т-диаграмма $X$ содержит более
четырех точек. Если из нее удалить черную точку $A_i$ и соответствующую ей белую $B$, то
получится вновь \textsl т-диаграмма.
\label{induct}
\end{lemma}

Доказательство. Воспользуемся предложением \ref{equal1enough}.
Треугольники, образованные черными точками в диаграмме $X'=X \setminus \{A_i, B\}$, были
треугольниками, образованными черными точками и в диаграмме $X$. Там они содержали
ровно одну белую точку, причем это была не точка $B$, которая соответствовала $A_i$
(см. замечание \ref{insideonly}). Значит и в каждом треугольнике, образованном черными точками $X'$,
содержится ровно одна белая точка диаграммы $X'$.\qed

\bigskip

Дадим два эквивалентных определения.

\begin{dfn}
Про белую точку т-диаграммы, которая одновременно лежит в граничных треугольниках
вершин $A_i$ и $A_{i+1}$, будем говорить, что она \textit{прилежит к стороне $A_iA_{i+1}$}.
\end{dfn}

\begin{dfn}
Будем говорить, что белая точка $B$ \textit{прилежит к стороне $A_iA_{i+1}$}, если она лежит
внутри всех треугольников вида $A_iA_{i+1}A_j$, и только в них. Это определение далее будем называть
альтернативным.
\end{dfn}

Подсчет числа граничных треугольников позволяет доказать следующую лемму.

\begin{lemma}
Для \textsl т-диаграммы справедливы следующие утверждения.
\begin{enumerate}
\item Если $d>0$, то существуют хотя бы 2 белые точки, прилежащие к сторонам, а при $d=0$ такая
точка ровно одна.
\item Если $d>0$, то белая точка не может соответствовать более, чем двум черным.
\item Белая точка соответствует двум данным черным тогда и только тогда, когда она прилежит к стороне
$(d+3)$-угольника, соединяющей эти черные точки, а сами они~--- соседние на границе
$(d+3)$-угольника.\qed
\end{enumerate}
\end{lemma}

\begin{dfn}
Две т-диаграммы $X$ и $Y$ назовем \textit{эквивалентными по диагоналям}, если  между
$X$ и $Y$ существует такая биекция $\varphi$, что
\begin{enumerate}
\item $\varphi$\label{equivcolor} сохраняет цвет точек,
\item $\varphi$\label{equivorder} сохраняет неориентированный циклический порядок черных точек,
\item белая\label{equivtriangs} точка $B$ лежит в выпуклой оболочке подряд идущих черных вершин $\{A_i\}$ в
диаграмме $X$ тогда и только тогда, когда $\varphi (B)$ лежит в выпуклой оболочке
$\{\varphi (A_i) \}$ в диаграмме $Y$.
\end{enumerate}
\end{dfn}

\begin{rem}
Условие \ref{equivtriangs} из этого определения  допускает следующую эквивалентную формулировку $3'$:
белая точка $B$ лежит в выпуклой оболочке трех черных вершин $\{A_i\}$ в  диаграмме $X$
тогда и только тогда, когда $\varphi (B)$ лежит в выпуклой оболочке $\{\varphi (A_i) \}$
в диаграмме $Y$.\label{rem17}
\end{rem}

В дальнейшем мы покажем, что эквивалентность двух Т-диаграмм по диагоналям равносильна
их комбинаторной эквивалентности.

\begin{prop}
Пусть в \textsl т-диаграмме $X$ белая точка $B$ прилежит к стороне $A_iA_{i+1}$.
Тогда \textsl т-диаграммы $X \setminus \{A_i, B \}$  и  $X \setminus \{A_{i+1}, B\}$ эквивалентны
по диагоналям.
\end{prop}
Доказательство. Пусть биекция $\varphi$ между диаграммами $X \setminus \{A_i, B \}$ и
$X \setminus \{A_{i+1}, B\}$ точке $A_{i+1}$ ставит в соответствие точку $A_i$, а
на остальных точках $\varphi$ --- тождественное отображение.

Свойства \ref{equivcolor} и \ref{equivorder} для $\varphi$ очевидно выполнены.
Для проверки свойства \ref{equivtriangs} воспользуемся замечанием \ref{rem17}.
Рассмотрим некоторую белую точку $C$. Для треугольников, никакая
вершина которых не совпадает ни с $A_i$, ни с $A_{i+1}$, очевидно выполнение свойства $3'$.
Пусть после удаления $A_i$ белая точка $C$ оказалась в треугольнике $\triangle A_{i+1}A_jA_k$. Покажем,
что если бы мы удалили не $A_i$, а $A_{i+1}$, она оказалась бы в $\triangle A_iA_jA_k$.
Действительно, рассмотрим четырехугольник $A_iA_{i+1}A_jA_k$ (не ограничивая общности, будем
считать, что $A_iA_j$~--- его диагональ). Внутри него лежат 2 белые точки. Одна из них~---
точка $C$, которая лежит внутри $\triangle A_{i+1}A_jA_k$, а другая~--- белая
точка $B$, прилежащая к стороне $A_iA_{i+1}$, которая по альтернативному определению лежит
внутри $\triangle A_iA_{i+1}A_k$. Значит, внутри $\triangle A_iA_jA_k$ лежит одна из этих
белых точек, причем  $B$ там лежать не может. Значит, там лежит точка $C$.\qed

\bigskip
Перейдем к рассмотрению характеристических деревьев т-диаграмм.

\begin{Ex}
На рис.~\ref{treesample} изображены 3-деревья, являющиеся характеристическими для диаграмм
Гейла, приведенных на рис.~\ref{FirstSamples}.
\end{Ex}

Следующая теорема позволяет параметризовать т-диаграммы (а значит, и Т-диаграммы) 3-деревьями.

\begin{thm}\mbox{}\label{characteristictree}
\begin{enumerate}
\item Для любой \textsl т-диаграммы
существует характеристическое дерево. Оно единственно с точностью до одновременной смены
циклического порядка во всех вершинах на противоположный.

\item Для любого 3-дерева с числом вершин $D+4$ существует \textsl т-диаграмма, для которой данное
дерево является характеристическим. Она единственна с
точностью до эквивалентности по диагоналям.
\end{enumerate}
\end{thm}

\begin{fig}
$$
\begin{array}{cc}
\myincludegraphics{szsfigures.11}\label{treesample} & \myincludegraphics{szsfigures.12}\\[10pt]
\myincludegraphics{szsfigures.13} & \myincludegraphics{szsfigures.14}
\end{array}
$$
\end{fig}

В рамках доказательства теоремы будем считать, что на черных точках выбран циклический
порядок (а не только неориентированный циклический порядок), т.~е. зафиксировано направление
\glqq по часовой стрелке\grqq{} на границе $(d+3)$-угольника. Кроме того, в определении понятия
\glqq быть характеристическим деревом\grqq{}
в свойстве 2 потребуем совпадения именно циклических порядков, а не только
неориентированных циклических порядков.

Сделаем вначале несколько вспомогательных замечаний.

\begin{lemma}
Пусть\label{WcorrBintree} дана \textsl т-диаграмма, и пусть некоторое 3-дерево является для нее
характеристическим.
Тогда следующие утверждения равносильны:
\begin{enumerate}
\item Белая точка $B$ соответствует черной точке $A_i$.
\item Висячая вершина дерева, соответствующая $A_i$, соединена ребром с вершиной, соответствующей
$B$.\qed
\end{enumerate}
\end{lemma}

\begin{cor}
Пусть\label{rogatkag} дана \textsl т-диаграмма, и пусть некоторое 3-дерево является для нее
характеристическим. Тогда следующие утверждения равносильны:
\begin{enumerate}
\item Белая точка $B$ прилежит к стороне $A_iA_{i+1}$.
\item Вершина дерева, соответствующая $B$, соединена ребрами с вершинами, соответствующими
$A_i$ и $A_{i+1}$.\qed
\end{enumerate}
\end{cor}

\begin{lemma}
\textbf{\textup{Лемма об индукции для характеристического дерева.\label{charactinduct}}}
Пусть черные точки $A_1, \ldots, A_{d+3}$ \textsl т-диаграммы $X$ занумерованы согласно обходу
$(d+3)$-угольника по часовой стрелке.
Пусть точка $B$ прилежит к стороне $A_iA_{i+1}$. Пусть $\varphi$~---
некоторая биекция между множеством $X$ и множеством вершин некоторого 3-дерева $G$.
Пусть вершина дерева $\varphi (B)$ степени 3 соединена ребрами с висячими
вершинами $\varphi (A_i)$ и $\varphi(A_{i+1})$, а циклический
порядок на ребрах, выходящих из $\varphi (B)$,
согласован с циклическим порядком на черных точках диаграммы
(т.~е.
правильно
задает направление обхода по часовой стрелке для стороны $A_iA_{i+1}$).

Рассмотрим новую \textsl т-диаграмму $X' = X  \setminus \{ A_i, B\}$ и новое дерево
$G'= G \setminus \{\varphi (A_i), \varphi (A_{i+1})\}$. Новую биекцию $\varphi '$ зададим
следующим образом.

Вершину $\varphi (B)$, которая в дереве $G'$ является висячей, поставим в соответствие
точке $A_{i+1}$. Все остальные вершины оставим поставленными в соответствие с помощью биекции
$\varphi$.

Тогда  дерево $G$ является характеристическим для исходной диаграммы
$X$ тогда и только тогда, когда новое дерево $G'$ является
характеристическим для новой \textsl т-диаграммы $X'$.
\end{lemma}

Доказательство состоит в аккуратной проверке выполнения условий из определения
характеристического дерева.\qed

Докажем первое утверждение теоремы (существование и единственность характеристического
дерева) индукцией по $d$. Если $d=0$, то единственно возможная
т-диаграмма~--- треугольник с белой точкой внутри. Единственно возможное 3-дерево из 4
вершин~--- дерево с одной вершиной степени 3, из которой идут 3 ребра в вершины
степени 1.

Пусть теперь $d>0$.

Пусть для т-диаграммы $X$ нашлось характеристическое дерево.
Для краткости будем обозначать его вершины и точки диаграммы одними и теми же буквами. Пусть белая
вершина $B$ прилежит к стороне $A_iA_{i+1}$. Тогда (см. следствие \ref{rogatkag}) из вершины дерева $B$
идут 2 ребра в висячие вершины $A_i$ и $A_{i+1}$. Удалим из дерева
эти висячие вершины, а из т-диаграммы  удалим точки $A_i$ и $B$.

По предположению индукции,
для полученной диаграммы характеристическое дерево существует и единственно,
причем, по лемме~\ref{charactinduct} об индукции для характеристического дерева это
и будет то новое дерево, которое мы получили.

Тем самым, характеристическое дерево для $X$ восстанавливается однозначно, но пока с
точностью до циклического порядка при вершине $B$. Этот циклический порядок тоже
восстанавливается однозначно, т.~к. вершины $A_i$ и $A_{i+1}$ должны идти по часовой
стрелке.

Существование напрямую следует из леммы об индукции для характеристического дерева.
Надо
только добавить к характеристическому дереву для диаграммы Гейла без 2 вершин 2 новые вершины
степени 1.
Первое утверждение теоремы доказано.

\begin{lemma}
Пусть $A_i$~--- черная точка в \textsl т-диаграмме $X$. Тогда к $X$ можно добавить черную точку
$A_{i+1/2}$ и белую точку $B$ так, что получится новая \textsl т-диаграмма, в которой
\begin{enumerate}
\item точка  $A_{i+1/2}$ лежит\label{l317near} между  $A_i$ и  $A_{i+1}$ на границе выпуклой оболочки множества
черных точек, и
\item точка $B$ прилежит к стороне $A_iA_{i+1/2}$.
\end{enumerate}
\label{inductbw}
\end{lemma}

Доказательство. Рассмотрим  прямую, проходящую через $A_i$ и оставляющую все остальные
точки диаграммы $X$  по одну сторону от себя.
Выберем на этой прямой достаточно близкую к $A_i$ новую черную точку $A_{i+1/2}$ такую,
чтобы было выполнено условие \ref{l317near} и чтобы каждая белая точка лежала бы
по ту же сторону от прямых, соединяющих вершины $(d+3)$-угольника с выбранной точкой, что и от
прямых, соединяющих эти вершины с $A_i$.
Добавим белую точку $B$ в пересечение треугольников $A_{i-1}A_iA_{i+1/2}$ и $A_iA_{i+1/2}A_{i+1}$.
Точка $B$ попадет в треугольники со стороной $A_iA_{i+1/2}$, и только в них.\qed

\begin{dfn}
Назовем вершину 3-дерева \textit{предвисячей}, если она соединена ребрами с 2 висячими
вершинами.
\end{dfn}

Докажем второе утверждение теоремы индукцией по $d$. База индукции ($d=0$) очевидна,
а леммы \ref{charactinduct} и \ref{inductbw} позволяют совершить индукционный переход.
Теорема доказана.\qed

\bigskip

Теорема \ref{characteristictree} позволяет вычислить количество Т-диаграмм для данного $d$. Выяснив,
как устроены автоморфизмы характеристических деревьев, можно получить следующую формулу, вывод
которой мы опускаем:

\begin{thm}
Число комбинаторно различных Т-диаграмм  для данного $d>0$ равно
$$
\frac{T_{d+1}}{2(d+3)}+\frac{3T_{(d+3)/2-1}}4+\frac{T_{(d+3)/3-1}}3+\frac{T_{d/2}}2.
$$
Здесь $T_x$ равно 0, если $x$~--- нецелое, и $x$-тому числу Каталана $C_{2x}^x/(x+1)$, если
$x$~--- целое.\qed
\end{thm}

\section{Единственность T-диаграммы Гейла для данного Т-многогранника}

\begin{dfn}
Пусть имеется некоторая луночка. Назовем вершины $A_i$ и $A_j$ ее \textit{крайними вершинами},
если все ее вершины~--- это все подряд идущие вершины $(d+3)$-угольника от $A_i$ до $A_j$
включительно, и только они.
\end{dfn}

Если луночка~--- не весь $(d+3)$-угольник, то ее крайние вершины определены однозначно.
Обратно, если имеются 2 несоседние вершины $A_i$ и $A_j$, то существуют ровно 2 луночки,
для которых $A_i$ и $A_j$ являются крайними вершинами.  Если же луночка~--- весь
$(d+3)$-угольник, то ее крайними вершинами могут служить любые 2 соседние вершины,
которые являются вершинами ее и еще одной (имеющей 2 вершины) луночки.

\textit{Не-гранью} многогранника $K$ будем называть набор его вершин, который не
является множеством вершин ни для какой грани многогранника $K$.

\begin{lemma}
\textbf{\textup{Критерий не-грани из $d+1$ вершины на языке Т-диаграммы.}} Пусть $X$~---
Т-диаграмма Т-многогранника $X^*$. Множество $M^*$ из $d+1$ вершины образует не-грань
Т-многогранника $X^*$ тогда и только тогда, когда множество $X \setminus M$ таково:
черные точки этого множества ~--- это вершины некоторой луночки $L$, а белые точки~---
это те и только те белые точки диаграммы  $X$, которые не лежат в $L$.\label{CoFace}
\end{lemma}

Доказательство. Пусть множество $X\setminus M$ указанного
вида. Тогда все белые вершины лежат внутри другой из 2 луночек с теми же крайними вершинами,
что $L$ 
(если $L$~--- весь $(d+3)$-угольник, то белых вершин в $X\setminus M$ нет вообще).
Следовательно, множество $X\setminus M$ не обладает свойством \bw.

Обратно, пусть некоторое $(d+1)$-элементное подмножество $M$ диаграммы $X$ порождает
не-грань. Рассмотрим точки множества $X\setminus M$. Пусть среди них $k$ черных.
Поскольку образуемый ими $k$-угольник $K$ можно разбить на $k-2$ треугольника, внутри
него находятся $k-2$ белые точки Т-диаграммы. Все они принадлежат $M$ (иначе было бы
выполнено свойство \bw{}). Следовательно, вне $k$-угольника лежат $d+1-(k-2)=d-k+3$
белых точек, ни одна из которых не принадлежит $M$ (иначе в $X \setminus M$ было
бы менее $d+3$ вершин). Допустим, что черные точки из  $X \setminus M$ не образуют луночку,
то есть идут не подряд. Тогда $(d+3)$-угольник (с вершинами в черных точках) при удалении из
него $K$ распадется на не менее чем две луночки, каждая из которых содержит белые точки из
$X\setminus M$ (см.~рис.~\ref{CoFaceFig}). Следовательно, множество $X\setminus M$ обладает
свойством \bw{}. Противоречие.\qed

\begin{fig}
$$
\myincludegraphics{szsfigures.15}\label{CoFaceFig}
$$
\end{fig}

\begin{dfn}
Назовем не-грань Т-многогранника из более чем $d+1$ вершины \textit{неособой}, если среди вершин,
которые в нее не входят, есть хотя бы две черных и хотя бы одна белая. Все остальные
не-грани из более чем $d+1$ вершины назовем \textit{особыми}.
\end{dfn}

\begin{lemma}
\textbf{\textup{Критерий\label{CoFaceMore} неособой не-грани из более чем $d+1$ вершины.}} 
Множество $M^*$ из более
чем $d+1$ вершины порождает неособую не-грань Т-многогранника $X^*$ тогда и только тогда, когда
выполнены следующие условия.
\begin{enumerate}
\item Множество $X\setminus M$ содержит хотя бы две черные точки.
\item Множество $X\setminus M$ содержит хотя бы одну белую точку.
\item При\label{CoFaceMore3} удалении  из $(d+3)$-угольника Т-диаграммы $X$ выпуклой оболочки множества черных точек
из $X\setminus M$, белые точки из $X\setminus M$  лежат только в одной из луночек, на которые
распадется $(d+3)$-угольник.
\end{enumerate}
\end{lemma}

Доказательство.
Условие \ref{CoFaceMore3} равносильно тому, что множество $X\setminus M$ не обладает свойством
\bw.\qed

\begin{lemma}
\textbf{\textup{Критерий особой не-грани из более чем $d+1$ вершины.}} Множество\label{CoFaceMoreGaiane}
$M^*$ из более
чем $d+1$ вершины порождает особую не-грань Т-многогранника $X^*$ в следующих и только следующих
случаях.
\begin{enumerate}
\item Множество $X \setminus M$ содержит не более одной черной точки.
\item В множестве $X \setminus M$ все точки черные.\qed
\end{enumerate}
\end{lemma}

\bigskip

Для точек Т-диаграммы $R_1,\ldots,R_r$ 
обозначим через
$\mathcal N(R_1,\ldots, R_r)$ 
число не-граней из $d+1$ вершины,
содержащих вершины $R_1^*,\ldots,R_r^*$.

Пусть $A_i\ne A_j$~--- черные точки некоторой Т-диаграммы $X$.
Введем следующие обозначения: пусть $l=l(A_i, A_j)$~--- число
вершин $(d+3)$-угольника, лежащих строго между $A_i$ и $A_j$, если идти от $A_i$ к $A_j$ против
часовой стрелки (или, как будем говорить, по \textit{левому пути}). Аналогично,  $r=r(A_i, A_j)$~---
число вершин $(d+3)$-угольника, лежащих строго между $A_i$ и $A_j$, если идти от $A_i$ к $A_j$ по
часовой стрелке (или по \textit{правому пути}).
Очевидно, $l+r=d+1$. Иногда мы будем использовать выражение \glqq число $l$ (или $r$) для вершины
$A_j$, считая от вершины $A_i$\grqq. В этих обозначениях верна следующая лемма.

\begin{lemma}
$\mathcal N (A_i, A_j)=C_l^2+C_r^2$. (Здесь и далее мы считаем $C_n^m=0$
при $m>n$ или $m<0$.)\label{countBB}
\end{lemma}

Доказательство. Действительно, рассмотрим произвольную не-грань $M^*$, содержащую вершины $A_i^*$ и $A_j^*$. При
удалении из диаграммы $X$ множества $M$, оставшиеся черные точки должны образовывать луночку
(см. критерий не-грани из $d+1$ вершины).
Т.~к. точки  $A_i$ и $A_j$ удалены, то точки этой луночки могут либо все лежать на
левом пути, либо все на правом. При этом луночка уже не может содержать все черные
точки диаграммы $X$, значит ее крайние вершины определены однозначно. А число способов выбрать
2 вершины на левом пути есть $C_l^2$, на правом $C_r^2$, всего получается $C_l^2+C_r^2$.
После того, как выбраны крайние вершины луночки (а значит и сама луночка, т.~к. $A_i$ не
может являться ее вершиной), белые точки
не-грани определены однозначно (по критерию не-грани из $d+1$ вершины).\qed

\bigskip
Пусть  $A_j\neq A_i$~--- черные точки Т-диаграммы. Для белой точки $B$, соответствующей точке
$A_j$ вычислим  число $\mathcal N (A_i, B)$ (т.~е. число не-граней, содержащих ребро $A_i^*B^*$).

Для этого введем некоторые обозначения.

\begin{fig}
$$
\myincludegraphics{szsfigures.16}
\label{adef}
$$
\end{fig}

Пусть точки $A_k$ и $A_{k+1}$ таковы, что $B$ лежит в (определенном однозначно)
треугольнике $\triangle A_jA_kA_{k+1}$. Определим число $a=a(A_i,A_j,B)$
следующим образом (см.~рис.~\ref{adef}).

Если одна из точек
$A_k$ и $A_{k+1}$ совпадает с $A_i$, то положим $a=0$. Иначе обе точки $A_k$ и $A_{k+1}$
лежат либо на левом пути, ведущем из $A_i$ к  $A_j$, либо на правом. Положим  число $a$
равным количеству вершин на этом пути между $A_i$
и той из вершин $A_k$ и $A_{k+1}$, которая встретится первой,
считая $A_k$ или $A_{k+1}$, но не считая $A_i$. В этих обозначениях верны
следующие две леммы:

\begin{lemma}
Точки $A_j$ и $A_k$ однозначно определяют следующий набор луночек, содержащих точку $B$.
\begin{enumerate}
\item Весь $(d+3)$-угольник;
\item луночка, крайние вершины которой не совпадают с $A_j$, но среди вершин которой
встречается $A_j$;
\item луночка, одна из крайних вершин которой есть $A_j$ и среди вершин которой содержатся
$A_k$ и $A_{k+1}$.\qed
\end{enumerate}
\label{KnowK}
\end{lemma}

\begin{lemma}
$\mathcal N (A_i, B)=lr+a$.\label{conutBW}
\end{lemma}

Доказательство. Рассмотрим произвольную не-грань $M^*$, содержащую вершины $A_i^*$ и $B^*$.
Содержать вершину $A_j^*$ она не может~--- иначе по лемме~\ref{induct} об индукции, множество
$X \setminus \{A_j, B\}$, состоящее из $2d+2$ точки, являлось бы Т-диаграммой. Если из него удалить еще $d-1$
точку, множество оставшихся точек будет обладать свойством \bw.

Итак, точка $A_j$ не принадлежит  $M$ и является вершиной некоторой луночки $L$, образованной
точками из $X \setminus M$. Эта луночка содержит точку $B$.

Поскольку луночка $L$ не совпадает со всем $(d+3)$-угольником, ее крайние вершины (скажем,
$A_m$ и $A_n$) определены однозначно.

Заметим, что если имеется произвольная луночка $L$, такая что среди ее вершин нет $A_i$, но она
содержит $B$, то, согласно критерию не-грани из $d+1$ вершины, эта луночка задает ровно одну
не-грань из $d+1$ вершины, содержащую ребро $A_i^*B^*$.

Посчитаем число возможных способов выбрать (неупорядоченную) пару вершин $(d+3)$-угольника
$A_m$ и $A_n$ так, чтобы они были крайними вершинами луночки $L$, обладающей указанными
свойствами. В соответствии с леммой \ref{KnowK}, рассмотрим 2 случая:

Случай 1. Ни одна из вершин $A_m$ и $A_n$ не совпадает с $A_j$. Поскольку $A_j$ является вершиной
$L$, а $A_i$~--- нет, то одна из точек $A_m$ и $A_n$ лежит
на левом пути, ведущем из $A_j$ к  $A_i$, а другая~--- на правом. Тогда число способов их выбрать
равно $lr$. Для каждой такой пары вершин существует ровно одна луночка, такая что $A_m$ и $A_n$
являются ее крайними вершинами, причем среди ее вершин нет $A_i$ (но есть $A_j$).

Случай 2. Одна из вершин $A_m$ и $A_n$, скажем, $A_m$, совпадает с $A_j$. Тогда,
согласно лемме \ref{KnowK}, среди вершин $L$ должны быть вершины $A_k$ и $A_{k+1}$. Значит,
вершина $A_n$ должна лежать на том же (левом или правом) пути от $A_i$ к $A_j$, что и
$A_k$ и $A_{k+1}$. Кроме того, если идти по этому пути от $A_i$ к $A_j$, то мы должны
сначала встретить вершину $A_n$, и только потом $A_k$ и $A_{k+1}$ (при этом $A_n$
может совпадать с той из вершин $A_k$ и $A_{k+1}$, которую мы встретим первой). По
определению, число таких вершин на (левом или правом) пути и есть $a$. Значит,
в этом случае существует ровно $a$ пар вершин $A_m$ и $A_n$, которые могут служить крайними
вершинами $L$. И снова для каждой из найденных пар вершин существует ровно одна луночка,
такая что $A_m$ и $A_n$ являются ее крайними вершинами, но среди ее вершин нет $A_i$.

Всего получили $lr+a$ возможных луночек $L$.\qed

\begin{prop}
Ребро Т-многогранника замечательное тогда и только тогда, когда на Т-диаграмме
одному из его концов соответствует черная точка, а другому~---  соответствующая
ей белая.\label{GoodCrit}
\end{prop}

Доказательство. Ребро, соединяющее 2 черные вершины, не может быть замечательным по лемме \ref{countBB}.
Действительно, поскольку $D\ge4$, то $l+r=d+1\ge3$, поэтому хотя бы одно из чисел $l$ и $r$
не меньше 2, значит $C_l^2+C_r^2>0$.

Ребро, соединяющее 2 белые вершины, не может быть замечательным, потому что все $d+1$ белые
точки образуют не-грань по критерию (лемма \ref{CoFace}).

Пусть теперь вершинами ребра служат черная точка $A_j$ и белая точка $C$, причем $C$ не
соответствует точке $A_j$. Пусть точке $A_j$  соответствует точка $B$.
Рассмотрим множество точек Т-диаграммы, состоящее из белой точки $B$ и всех черных
точек, кроме $A_j$. Для этого множества не выполнено условие \bw,
следовательно его дополнение $M$ порождает не-грань. Противоречие
замечательности.

В обратную сторону утверждение следует из леммы об индукции, как это уже было проверено
при доказательстве предыдущей леммы.\qed

Теперь все готово для доказательства следующей теоремы.

\begin{thm}\mbox{}
\begin{enumerate}
\item Комбинаторика Т-многогранника однозначно восстанавливается по классу эквивалентности по
диагоналям его Т-диаграммы Гейла.\label{equivTpolyTdiag}

\item Для любого Т-многогранника $X$ комбинаторика (в смысле эквивалентности по диагоналям)
его Т-диаграммы Гейла $X^*$ \textbf{однозначно} восстанавливается по
комбинаторным свойствам $X$.
\end{enumerate}
\end{thm}

Первая часть теоремы следует из критериев не-грани (см. леммы \ref{CoFace}, \ref{CoFaceMore} и 
\ref{CoFaceMoreGaiane}).

Докажем вторую часть теоремы.

Пусть имеется Т-многогранник $X^*$ и его Т-диаграмма $X$. Покажем, что цвета точек
$X$ однозначно определены  комбинаторикой $X^*$.

Благодаря предл. \ref{GoodCrit}, комбинаторика $X^*$ позволяет различать пары из
черной вершины и соответствующей ей белой. Пока мы не можем отличить черную вершину от
белой, если эта белая вершина соответствует ровно одной черной. Но если белая точка
соответствует сразу двум черным (такие белые точки существуют),
то из нее выходит 2 замечательных ребра, и мы можем отличить ее от
2 черных, которым она соответствует. Без ограничения общности будем считать, что это
вершины $A_1$ и $A_{d+3}$.
Зафиксируем  черную точку $C$ (отличную от $A_1$ и $A_{d+3}$) и
соответствующую ей белую точку $B$.

В этих обозначениях сформулируем и докажем ряд вспомогательных лемм.

\begin{lemma}
Число $\mathcal N(A_1, C)-\mathcal N(A_{d+3}, C)$ имеет другую четность\label{odddiff},
чем число $d+1$.
\end{lemma}
Доказательство.
Воспользуемся леммой \ref{countBB} и заметим, что при замене $A_1$ на
$A_{d+3}$ числа $l$ и $r$  меняются на 1.

Подробнее, пусть $l_1$ и $r_1$~--- числа $l$ и $r$ для вершины $C$, считая от
одной из вершин $A_1$ или $A_{d+3}$, а $l_2$ и $r_2$~--- числа
$l$ и $r$ для вершины $C$, считая от другой из них. Не ограничивая общности, будем считать, что
$l_1-l_2=1$, $r_1-r_2=-1$ (иначе поменяем местами числа $l_1$ и $l_2$, а так же
$r_1$ и $r_2$).
$|\mathcal N(A_1, C)-\mathcal N(A_{d+3}, C)|=|C_{l_1}^2+C_{r_1}^2-C_{l_2}^2-C_{r_2}^2|=
|C_{l_2+1}^2-C_{l_2}^2+C_{r_2-1}^2-C_{r_2}^2|=|C_{l_2}^1-C_{r_2-1}^1|=|l_2-r_2+1|$.
Т.~е. число $\mathcal N(A_1, C)-\mathcal N(A_{d+3}, C)$ имеет другую
четность, чем $l_2-r_2$, а значит и другую четность, чем $l_2+r_2=d+1$.\qed

\begin{lemma}
Число $\mathcal N(A_1, B)-\mathcal N(A_{d+3}, B)$ совпадает по четности\label{evendiff} с числом $d+1$.
\end{lemma}
Доказательство. Пусть $l_1$, $l_2$, $r_1$ и $r_2$~---
числа, описанные выше (причем $l_1-l_2=1$, $r_1-r_2=-1$). Положим
$a_1=a(A_1,C,B)$ и $a_2=(A_{d+3},C,B)$.

Числа $a_1$ и $a_2$ не могут быть равны 0 одновременно, т.~к. иначе вершина $B$ лежала бы в
$\triangle A_1A_{d+3}C$, в котором уже есть одна белая вершина, а именно прилежащая к
стороне $A_1A_{d+3}$. Поэтому число $a_1-a_2$ равно либо 1, либо $-1$.
Теперь утверждение следует из леммы \ref{conutBW}.\qed

Благодаря леммам \ref{odddiff} и \ref{evendiff}, комбинаторика Т-многогранника
однозначно определяет цвет каждого из концов любого замечательного
ребра.

Это значит, что цвета всех точек также определены.

Следовательно, для векторной диаграммы Гейла (определяемой однозначно) $\overline X$ про
каждый вектор известно, какого цвета точку он порождает на аффинной Т-диаграмме Гейла.

Следовательно, выбор плоскости $e$ (см.~рис.~\ref{vecaffdiag}) в большой степени предопределен
и не влияет на комбинаторику Т-диаграммы.\qed

\begin{cor}
Следующие 2 утверждения равносильны:
\begin{enumerate}
\item Две Т-диаграммы  $X$ и $Y$ эквивалентны по диагоналям.
\item Две Т-диаграммы  $X$ и $Y$ комбинаторно эквивалентны.
\end{enumerate}
\label{combequiv}
\end{cor}
Доказательство. Следствие $2\Rightarrow1$ очевидно, докажем $1\Rightarrow2$.

Пусть две Т-диаграммы $X$ и $Y$ эквивалентны по диагоналям. По первой части теоремы, отсюда
следует, что $X^*$ и $Y^*$ комбинаторно эквивалентны как многогранники. Но они смежностные,
поэтому, как известно (см. вступление), их множества вершин комбинаторно эквивалентны как
множества точек. Значит, соответствующие векторные диаграммы Гейла $\overline X$ и
$\overline Y$ комбинаторно эквивалентны.
Кроме того, мы доказали, что если $X^*$ и $Y^*$~--- Т-многогранники, то про каждую их
вершину известно, какого цвета точку она даст в Т-диаграмме. Т.~к. $\varphi$ дает комбинаторный
изоморфизм $X^*$ и $Y^*$ как множеств точек (и, тем более, как многогранников), то $\varphi$
сохраняет цвет той точки, которая получится из вектора в $\overline X$ и в $\overline Y$.
Поэтому $X$ и $Y$ комбинаторно эквивалентны.\qed
\begin{cor}
Группы комбинаторных  автоморфизмов $X$ и $X^*$ канонически изоморфны.\qed
\end{cor}

\begin{cor}
Группа автоморфизмов Т-многогранника не действует на его вершинах транзитивно (в отличие
от циклического многогранника).\qed
\end{cor}

\begin{cor}\label{TpolyTree}\mbox{}
\begin{enumerate}
\item Любому Т-многограннику соответствует 3-дерево (характеристическое дерево для его Т-диаграммы
Гейла). Оно единственно с точностью до смены циклического порядка
вершин на противоположный.
\item Любое 3-дерево с не менее, чем 8 вершинами, задает единственный с точностью
до комбинаторной эквивалентности Т-многогранник.\qed
\end{enumerate}
\end{cor}

Для восстановления Т-диаграммы Гейла по комбинаторике Т-многогранника мы использовали
достаточно сложные инварианты (число не-граней, содержащих данное ребро).
Следующая теорема (мы ее приводим без доказательства) утверждает, что более простым
инвариантом обойтись не удастся.

\begin{thm}
В любом Т-многограннике размерности $D$ число граней с
$m$ вершинами, содержащих данную вершину, зависит только от $D$ и $m$,
но не от диаграммы и выбранной вершины. Оно совпадает с соответствующим числом
для циклического многогранника.\label{EqualIncedence}
\end{thm}

Доказательство основано на леммах \ref{CoFace}, \ref{CoFaceMore} и \ref{CoFaceMoreGaiane}.\qed

\bigskip

\footnotesize \textit{E-mail:} \texttt{deviatov1@rambler.ru}
\end{document}